\documentclass[12pt,twoside]{article}%
\usepackage{enumerate}
\usepackage{enumitem}
\usepackage{amsmath,amssymb,amsfonts}
\usepackage{amsthm}%
\usepackage{amsmath}%
\usepackage{amsfonts}%
\usepackage{amssymb}%
\usepackage{graphicx}
\usepackage{color}
\providecommand{\U}[1]{\protect\rule{.1in}{.1in}}
\newtheorem{thm}{Theorem}[section]
\newtheorem{lem}[thm]{Lemma}

\newtheorem{prop}[thm]{Proposition}
\newtheorem{corollary}[thm]{Corollary}
\newtheorem{claim}[thm]{Claim}

\newenvironment{pf}[1][\bfseries Proof]{\noindent{#1.} }{\hfill \rule{0.5em}{0.5em}\\}
\newcommand{\fim}{\hfill\rule{2mm}{2mm}}
\numberwithin{equation}{section}

\begin{document}

\title{Multiple solutions for a class of
quasilinear problems involving variable exponents \thanks{Partially supported by INCT-MAT and
PROCAD}}
\author{Claudianor O. Alves\thanks{C.O. Alves was partially supported by CNPq/Brazil 303080/2009-4,
e-mail:coalves@dme.ufcg.edu.br} \,\,\, and \,\,\, Jos\'{e} L. P.
Barreiro\thanks{e-mail:lindomberg@dme.ufcg.edu.br}\\Universidade Federal de Campina Grande\\Unidade Acad\^emica de Matem\'atica \\CEP:58429-900, Campina Grande - PB, Brazil.}
\date{}
\maketitle

\begin{abstract}
In this paper we prove the multiplicity of solutions for a class of quasilinear problems in $ \mathbb{R}^{N} $ involving variable exponents. The main tool used is in the proof are the direct methods, Ekeland's variational principle and some properties related to Nehari manifold.
\end{abstract}



\section{Introduction}

In this paper, we consider the existence and multiplicity of solutions for the following class of quasilinear problems involving variable exponents
\begin{align}
\left\{
\begin{array}
[c]{rcl}%
-\Delta_{p(x)} u + \vert u \vert^{p(x) - 2} u & = & \lambda g(k^{-1}x)
\vert u \vert^{q(x) - 2}u + f(k^{-1}x) \vert u
\vert^{r(x) - 2}u , \quad
\mathbb{R}^{N}\\
\mbox{}\\
u \in W^{1, p(x)}(\mathbb{R}^{N}) &  &
\end{array}
\right. \tag{$ P_{\lambda, k} $}\label{Plkm}%
\end{align}
where $\lambda $ and $k$ are nonnegative parameters with $k \in \mathbb{N}$, the operator $\Delta_{p(x)}u = \mathrm{div}\left( |\nabla u|^{p(x) - 2}\nabla u \right) $ is called $p(x)$-Laplacian, which is a natural extension of the $p$-Laplace operator, with $p$ being a positive constant. We assume that $p,q, r : \mathbb{R}^{N}
\to \mathbb{R} $ are positive Lipschitz continuous functions, which are $ \mathbb{Z}^{N}$-periodic and verify
\begin{align}
1 <
p_{-} \leq p(x) \leq p_{+} < q_{-} \leq q(x) \leq r(x) \ll p^{*}(x) \mbox{ a.e. in }
\mathbb{R}^{N}, \tag{$p_{1}$} \label{p1}
\end{align}
where $p_{+} =  \mbox{ess}\sup_{x \in \mathbb{R}^{N}}{p(x)}$, $p_{-} =  \mbox{ess}\inf_{x \in \mathbb{R}^{N}}{p(x)}$ and
$$
p^{*}(x) = \left\{
\begin{array}{l}
Np(x) / (N - p(x)) ~~ ~\mbox{if} ~~  p(x) < N  \\
\mbox{}\\
+\infty ~~ \mbox{if} ~~ p(x) \geq N .
\end{array}
\right.
$$
Moreover, we say that a measurable function $h: \mathbb{R}^{N} \to \mathbb{R}$ is $ \mathbb{Z}^{N}$-periodic if
$$
h(x + z) = h(x) \,\,\, \forall x \in \mathbb{R}^{N} \,\,\, \mbox{and} \,\,\, \forall z \in \mathbb{Z}^{N},
$$
and the notation $u \ll v$ means that $\displaystyle \inf_{x \in \mathbb{R}^{N}} (v(x) - u(x)) > 0$.

\medskip

Related to functions $f$ and $g $, we suppose that they are nonnegative continuous functions verifying the following conditions:

\begin{enumerate}[label={($H\arabic{*}$})]
\setcounter{enumi}{0}
\item\label{H2}
\begin{align*}
\lim_{\vert x \vert\rightarrow\infty} {g(k^{-1}x)} = 0;
\end{align*}

\item\label{H3}
There exist $\ell$ points $a_{1}, a_{2}, \cdots,a_{\ell} $ in
$\mathbb{Z}^{N} $ with $a_1=0$ such that
\[
1 = f(a_{i}) = \max_{\mathbb{R}^{N}}f(x), \text{ for }1 \leq i \leq\ell.
\]
Moreover, $0 < f_{\infty} < f(x) $ for any $x \in\mathbb{R}^{N} $ and
\[
\lim_{\vert x \vert\rightarrow\infty} f(x) = f_{\infty} .
\]

\end{enumerate}

The variable exponents problems appear in a lot of applications, the reader can find in R\r{u}\v{z}i\v{c}ka \cite{Ru} and Krist\'aly, Radulescu \& Varga in \cite{KRV} several models in mathematical physics where this class of problem appear. In recent years, such problems have attracted an increasing attention, we would like to mention \cite{alves08, AlvesFerreira1, AlvesSouto,chabrowki, FanHan, MR}, and also the survey papers \cite{AS,DHN,S} for the advances and references in this area.

\medskip

The problem $(P_{\lambda,k})$ has been considered in the literature for the case where the exponents are constants, see for example, Adachi \& Tanaka \cite{AT}, Cao \& Noussair \cite{CN}, Cao \& Zhou \cite{Cao}, Hirano \cite{H1}, Hirano \& Shioji \cite{HS}, Hu \& Tang \cite{HuTang}, Jeanjean \cite{jeanjean}, Lin \cite{Lin12}, Hsu, Lin \& Hu \cite{Hsu1}, Tarantello \cite{T}, Wu \cite{Wu1, Wu2} and their references.

\medskip

In Cao \& Noussair \cite{CN}, the authors have studied the existence and multiplicity of positive and nodal solutions for the following problem
$$
\left\{
\begin{array}{l}
-\Delta u + u  =  f(\epsilon x) \vert u \vert^{r - 2}u  \,\,\, \mbox{in} \,\,\, \mathbb{R}^{N}\\
\mbox{}\\
u \in H^{1,2}(\mathbb{R}^{N}),
\end{array}
\right. \eqno{(P_{1})}
$$
where $\epsilon$ is a positive real parameter, $r \in (2,2^{*})$  and $f$ verifies condition $(H2)$. By using variational methods, the authors showed the existence of at least $\ell$ positive solutions and $\ell$ nodal solutions if $\epsilon$ is small enough. After, Wu in \cite{Wu1} considered the perturbed problem
$$
\left\{
\begin{array}{l}
-\Delta u + u  =  f(\epsilon x) \vert u \vert^{r - 2}u + \lambda g(\epsilon x)|u|^{q-2}u \,\,\, \mbox{in} \,\,\, \mathbb{R}^{N}\\
\mbox{}\\
u \in H^{1,2}(\mathbb{R}^{N}),
\end{array}
\right. \eqno{(P_{2})}
$$
where $\lambda$ is a positive parameter and $q \in (0,1)$. In \cite{Wu1}, the authors showed the existence of at least $\ell$ positive solutions for $(P_2)$ when  $\epsilon$ and $\lambda$ are small enough.

\medskip

In Hsu, Lin \& Hu \cite{Hsu1}, the authors have considered the following class of quasilinear problems
$$
\left\{
\begin{array}{l}
-\Delta_p u + |u|^{p-2}u  =  f(\epsilon x) \vert u \vert^{r - 2}u + \lambda g(\epsilon x) \,\,\, \mbox{in} \,\,\, \mathbb{R}^{N}\\
\mbox{}\\
u \in W^{1,p}(\mathbb{R}^{N})
\end{array}
\right. \eqno{(P_{3})}
$$
with $N \geq 3$ and $2 \leq p < N$. In that paper, the authors have proved the same type of results found in \cite{CN} and \cite{Wu1}.

\medskip

Motivated by results proved in \cite{CN}, \cite{Hsu1} and \cite{Wu1}, we intend in the present paper to prove the existence of multiple solutions for problem (\ref{Plkm}), by using the same type of approach explored in that papers. However, once that we are working with variable exponents, some estimates that hold for the constant case are not immediate for the variable case, and so, a careful analysis is necessary to get some estimates. Here, for example, we were able to prove our results by assuming that some exponents are periodic and $k \in \mathbb{N}$.

\vspace{0.5 cm}

Our main result is the following

\begin{thm}\label{T1}
Assume that (\ref{p1}) and \ref{H2}--\ref{H3} are satisfied. Then
there are $ \Lambda^{*} >0 $ and $k^* \in \mathbb{N}$ such that problem (\ref{Plkm}) admits at least $ \ell $ solutions for $ 0 \leq  \lambda < \Lambda^{*} $ and $k \geq k^{*}$.
\end{thm}

This paper is organized in the following way: In Section~\ref{expoents_variaveis}, we collect some preliminaries on variable exponent spaces that will be used throughout the paper, which can be found in \cite{AlvesFerreira1}, \cite{AlvesFerreira2}, \cite{chabrowki}, \cite{peter}  and  \cite{Fan2001a}. In Section 3, we show some technical results,  and finally in Section 5 we prove Theorem \ref{T1}.

\vspace{0.5 cm}
\noindent\textbf{Notation:} The following notation will be used in the present work:
\begin{itemize}
  \item $ C $ and $ c_{i} $ denote generic positive constants, which may vary from line to line.
  \item  We denote by $\int f $ the integral $\int_{\mathbb{R}^{N}}fdx$, for any measurable function $f$.
  \item $ B_{R}(z) $ denotes the open ball with center at $ z $ and radius $ R $ in $\mathbb{R}^{N}$.
	\item If $h$ is a bounded mensurable function, we denote by $h_+$ and $h_-$ the ensuing real numbers
	$$
	h_+=\mbox{ess}\sup_{x \in \mathbb{R}^{N}}{h(x)} ~~~~~~ \mbox{and} ~~~~~~ h_-=\mbox{ess}\inf_{x \in \mathbb{R}^{N}}{h(x)}.
	$$
	Moreover, we also denote by $h'(x)$ the conjugate exponent of $h(x)$ given by $h'(x)=\frac{h(x)}{h(x)-1}$.

\end{itemize}

\section{Preliminaries on Lebesgue and Sobolev spaces with variable exponent in $\mathbb{R}^{N}$}

\label{expoents_variaveis} In this section, we recall the definitions and some
results involving the spaces $L^{h(x)}(\mathbb{R}^{N}) $ and $W^{1,h(x)}%
(\mathbb{R}^{N}) $. We refer to \cite{peter,Fan2001a, Fan2001b, kovacik91} for the fundamental properties of these spaces.

Hereafter, let us denote by $L_{+}^{\infty}(\mathbb{R}^{N})$ the set
\[
L_{+}^{\infty}(\mathbb{R}^{N}) = \left\{  u \in L^{\infty}(\mathbb{R}^{N}) : \mbox{ess}\inf_{x
\in\mathbb{R}^{N}}u \geq1\right\}
\]
and we will assume that $h \in L_{+}^{\infty}(\mathbb{R}^{N})$.

The variable exponent Lebesgue space $L^{h(x)}(\mathbb{R}^{N}) $ is defined by
\[
L^{h(x)}(\mathbb{R}^{N}) = \left\{  u :\mathbb{R}^{N} \to \mathbb{R} \text{ is measurable }\, :
\, \, \,\, \int \vert u(x)\vert^{h(x)} < + \infty\right\},
\]
which is endowed with the norm
\[
\Vert u \Vert_{h(x)} = \inf\left\{  t > 0 : \int
\left\vert \frac{u(x)}{t}\right\vert ^{h(x)} \leq1\right\} .
\]
On space $L^{h(x)}(\mathbb{R}^{N})$, we consider the \textit{modular function} $\rho:
L^{h(x)}(\mathbb{R}^{N}) \to\mathbb{R}$ given by
\[
\rho(u) = \int |u(x)|^{h(x)} .
\]

\begin{prop}
\label{modular} Let $u \in L^{h(x)}(\mathbb{R}^{N}) $ and $\{u_{n}\}_{n \in\mathbb{N}}
\subset L^{h(x)}(\mathbb{R}^{N}) $. Then,

\begin{enumerate}
\item If $u \neq0 $, $\Vert u \Vert_{h(x)} = a
\Leftrightarrow\rho\left(  \frac{u}{a}\right)  = 1$.

\item $\Vert u \Vert_{h(x)} < 1 \quad(=1; > 1) \Leftrightarrow
\rho(u) < 1 (= 1; > 1) $;

\item $\Vert u \Vert_{h(x)} > 1 \Rightarrow\Vert u \Vert_{h(x)}^{h_{-}} \leq\rho(u) \leq\Vert u \Vert_{h(x)}^{h_{+}} $.

\item $\Vert u \Vert_{h(x)} < 1 \Rightarrow\Vert u \Vert_{h(x)}^{h_{+}} \leq\rho(u) \leq\Vert u \Vert_{h(x)}^{h_{-}} $.

\item $\displaystyle \lim_{n \to+\infty} \Vert u_{n} \Vert_{h(x)}
= 0 \Leftrightarrow\lim_{n \to\infty}\rho(u_{n}) = 0 .$

\item $\displaystyle \lim_{n \to+\infty} \Vert u_{n} \Vert_{h(x)}
= + \infty\Leftrightarrow\lim_{n \to\infty}\rho(u_{n}) = + \infty$.
\end{enumerate}
\end{prop}


\medskip
We have the following H\"{o}lder inequality for Lebesgue spaces with variable exponents.

\begin{prop}
[H\"{o}lder-type Inequality]Let $u \in L^{h(x)}(\mathbb{R}^{N})$ and $v \in L^{h^{\prime
}(x)}(\mathbb{R}^{N}) $. Then, $uv \in L^{1}(\mathbb{R}^{N})$ and
\begin{align*}
\int \vert u(x)v(x)\vert  \leq\left(  \frac{1}{h_{-}} + \frac
{1}{h^{\prime}_{-}} \right)  \Vert u \Vert_{h(x)} \Vert v
\Vert_{h^{\prime}(x)}.
\end{align*}
\end{prop}




%

The next three results are important tools  to study the properties of some energy functionals, and their
proofs can be found in \cite{AlvesFerreira2}.

\begin{lem} [Brezis-Lieb's lemma, first version] \label{Brezis-Lieb-1}
   Let $ ( \eta_n ) \subset L^{ h(x) } ( \mathbb R^N, \mathbb R^m ) $  with $ m \in \mathbb{N} $ verifying

   \begin{enumerate}
      \item [\emph{(i)}] $ \eta_n(x) \to \eta(x), \ \text{a.e. in} \ \mathbb R^N $;
      \item [\emph{(ii)}] $ \displaystyle \sup_{n \in \mathbb N } | \eta_n |_{ L^{ h(x) } ( \mathbb R^N, \mathbb R^m ) } < \infty $. \\
   \end{enumerate}
   Then, $ \eta \in L^{ h(x) } ( \mathbb R^N, \mathbb R^m ) $ and
   \begin{equation*}
      \int \left( \left| \eta_n \right|^{ h(x) } - \left| \eta_n - \eta \right|^{ h(x) } - \left| \eta \right|^{ h(x) } \right) \,= o_n(1).
   \end{equation*}
\end{lem}

\begin{lem} [Brezis-Lieb's lemma, second version] \label{Brezis-Lieb-2}
   Let $( \eta_n ) \subset L^{ h(x) } ( \mathbb R^N, \mathbb R^m ) $ with $ m \in \mathbb{N} $ verifying

   \begin{enumerate}
      \item [\emph{(i)}] $ \eta_n(x) \to \eta(x), \ \text{a.e. in} \ \mathbb R^N $;
      \item [\emph{(ii)}] $ \displaystyle \sup_{n \in \mathbb N } | \eta_n |_{ L^{ h(x) } ( \mathbb R^N, \mathbb R^m ) } < \infty $. \\
   \end{enumerate}
   Then
   \begin{equation*}
      \eta_n \rightharpoonup \eta \ \text{in} \ L^{ h(x) } ( \mathbb R^N, \mathbb R^m ).
   \end{equation*}
\end{lem}

The next proposition is a Brezis-Lieb type result.

\begin{lem} [Brezis-Lieb lemma, third version]
   Let   $ ( \eta_n ) \subset L^{ h(x) }  ( \mathbb R^N, \mathbb R^m ) $ with $ m \in \mathbb{N} $ such that

   \begin{enumerate} \label{Brezis-Lieb-3}
      \item [\emph{(i)}] $ \eta_n(x) \to \eta(x), \ \text{a.e. in} \ \mathbb R^N $;
      \item [\emph{(ii)}] $ \displaystyle \sup_{n \in \mathbb N } | \eta_n |_{ L^{ h(x) } ( \mathbb R^N, \mathbb R^m ) } < \infty $. \\
   \end{enumerate}
   Then
   \begin{equation*}
      \int \left| \left| \eta_n \right|^{ h(x)-2 } \eta_n - \left| \eta_n - \eta \right|^{ h(x)-2 } \left( \eta_n - \eta \right) - \left| \eta \right|^{ h(x)-2 } \eta \right|^{ h'(x) } \,  = o_n(1).
   \end{equation*}
 \end{lem}

The variable exponent Sobolev space $W^{1,h(x)}(\mathbb{R}^{N}) $ is defined by
\begin{align*}
W^{1,h(x)}(\mathbb{R}^{N}) = \left\{  u \in W^{1,1}_{loc}(\mathbb{R}^{N}) : u \in
L^{h(x)}(\mathbb{R}^{N}) \quad\text{ and } \quad| \nabla u | \in L^{h(x)}(\mathbb{R}^{N})
\right\}.
\end{align*}
The corresponding norm for this space is
\begin{align*}
\Vert u \Vert_{W^{1,h(x)}(\mathbb{R}^{N})} = \Vert u \Vert_{h(x)} +
\Vert\nabla u \Vert_{h(x)}.%
\end{align*}
The spaces $L^{h(x)}(\mathbb{R}^{N})$ and $W^{1,h(x)}(\mathbb{R}^{N})$ are separable and reflexive Banach spaces when $h_{-} >1$.

On space $W^{1,h(x)}(\mathbb{R}^{N})$, we consider the {\it modular function}  \linebreak $\rho_{1}: W^{1,h(x)}(\mathbb{R}^{N}) \to \mathbb{R}$ given by
\[
\rho_{1}(u) = \int \left(  |\nabla u(x)| ^{h(x)} + |u(x)| ^{h(x)}\right).
\]
If, we define
\begin{align*}
\Vert u \Vert = \inf\left\{ t > 0 : \int \frac{(|\nabla u |^{h(x)}+| u |^{h(x)})}{t^{h(x)} }  \leq 1\right\},
\end{align*}
then $ \Vert \cdot \Vert_{W^{1,h(x)}(\mathbb{R}^{N})} $ and $ \Vert \cdot \Vert $ are equivalent norms on $ W^{1,h(x)}(\mathbb{R}^{N}) $.
\begin{prop}\label{modular2}
Let $ u \in W^{1,h(x)}(\mathbb{R}^{N}) $ and $ \{u_{n}\} \subset W^{1,h(x)}(\mathbb{R}^{N}) $. Then, the same conclusion of Proposition \ref{modular} occurs replacing $\|\,\,\, \|_{h(x)}$ and $\rho$ by $\|\,\,\, \|$ and $\rho_1$ respectively.

\end{prop}
%
%



\section{Technical lemmas}
Associated with problem~(\ref{Plkm}), we have the energy functional
\linebreak$J_{\lambda, k}:W^{1,p(x)}(\mathbb{R}^{N}) \to \mathbb{R}$ defined
by
\[
J_{\lambda, k}(u) =\int\frac{1}{p(x)}\left(  |\nabla
u|^{p(x)}+|u|^{p(x)}\right)  -\lambda\int \frac{g(k^{-1}x)}%
{q(x)}|u|^{q(x)} - \int \frac{f(k^{-1}x)
}{r(x)}|u|^{r(x)} .
\]
It is easy to see that $J_{\lambda, k}\in C^{1}\left(  W^{1,p(x)}%
(\mathbb{R}^{N}),\mathbb{R}\right)  $ with
\begin{align*}
J'_{\lambda, k}(u) v  & =\int \left(  |\nabla
u|^{p(x)-2}\nabla u\nabla v+|u|^{p(x)-2}uv\right)  -\lambda\int g(k^{-1}x)|u|^{q(x)-2}uv\\
& \quad -\int f(k^{-1} x) |u|^{r(x)-2}uv  ,
\end{align*}
for any $u,v\in W^{1,p(x)}(\mathbb{R}^{N})$. Thus, the critical points of $ J_{\lambda, k} $ are (weak) solutions  of  (\ref{Plkm}). Since the functional $J_{\lambda, k}$ is not bounded from below on $ W^{1,p(x)}(\mathbb{R}^{N})$ , we will work on \emph{Nehari manifold} $ \mathcal{M}_{\lambda, k}$ associated with the functional $J_{\lambda, k}$, given by
$$
\mathcal{M}_{\lambda, k} = \left\{ u \in W^{1,p(x)}(\mathbb{R}^{N})\setminus\{ 0\}: J'_{\lambda, k}(u) u = 0 \right\}
$$
and with level
\[
c_{\lambda, k} = \inf_{u \in \mathcal{M}_{\lambda, k}} J_{\lambda, k}(u).
\]

Using well known arguments found in Willem \cite{willem}, it follows that $c_{\lambda, k} $ is the mountain pass level of functional $J_{\lambda, k}$.

For $ f \equiv 1 $ and $ \lambda = 0 $, we consider the problem
\begin{align}
\left\{
\begin{array}
[c]{rcl}%
-\Delta_{p(x)} u + \vert u \vert^{p(x) - 2} u & = & \vert u
\vert^{r(x) - 2}u , \quad
\mathbb{R}^{N}\\
\mbox{}\\
u \in W^{1, p(x)}(\mathbb{R}^{N}). &  &
\end{array}
\right. \tag{$ P_{\infty} $}\label{Poo}%
\end{align}
Associated with the problem~(\ref{Poo}), we have the energy functional \linebreak $ J_{\infty}:W^{1, p(x)}(\mathbb{R}^{N}) \to \mathbb{R}$ given by
\begin{align*}
J_{\infty}(u)   =\int \frac{1}{p(x)}\left(  |\nabla
u|^{p(x)}+|u|^{p(x)}\right)    -  \int \frac{1
}{r(x)}|u|^{r(x)},
\end{align*}
the level
$$
c_{\infty} = \inf_{u \in \mathcal{M}_{\infty}} J_{\infty}(u),
$$
and the Nehari manifold
$$
\mathcal{M}_{\infty} = \left\{ u \in W^{1,p(x)}(\mathbb{R}^{N})\setminus\{ 0\}: J'_{\infty}(u)u = 0 \right\}.
$$

For $ f \equiv f_{\infty} $ and $ \lambda = 0 $, we fix the problem
\begin{align}
\left\{
\begin{array}
[c]{rcl}%
-\Delta_{p(x)} u + |u|^{p(x) - 2} u & = & f_{\infty} |u|^{r(x) - 2}u, \quad
\mathbb{R}^{N}\\
\mbox{}\\
u \in W^{1, p(x)}(\mathbb{R}^{N}), &  &
\end{array}
\right. \tag{$ P_{f_{\infty}} $}\label{Pf00}%
\end{align}
and as above, we denote by $ J_{f_\infty}, c_{f_{\infty}}$ and $\mathcal{M}_{f_{\infty}}$ the energy functional, the mountain pass level and Nehari manifold associated with $ (P_{f_{\infty}}) $ respectively.

\vspace{0.5 cm}
The following result concerns the behavior of $J_{\lambda, k} $ on $ \mathcal{M}_{\lambda, k} $.

\begin{lem}
The functional $ J_{\lambda, k} $ is bounded from below on $ \mathcal{M}_{\lambda, k} $. Moreover, $ J_{\lambda, k} $ is coercive on $ \mathcal{M}_{\lambda, k} $.
\end{lem}
\noindent \textbf{Proof.}
For each  $ u \in \mathcal{M}_{\lambda, k} $,  $ J'_{\lambda, k}(u)  u = 0 $. Hence,
\begin{align*}
\lambda \int g(k^{-1}x) |u|^{q(x)}  = \int \left(  |\nabla
u|^{p(x)}+|u|^{p(x)}\right)   - \int f(k^{-1} x) |u|^{r(x)}.
\end{align*}
Note that
\begin{align*}
J_{\lambda, k}(u) &\geq  \frac{1}{p_{+}} \int \left(  |\nabla u|^{p(x)}+|u|^{p(x)}\right)   - \frac{\lambda}{q_{-}} \int  g(k^{-1}x) |u|^{q(x)} - \frac{1}{{r_{-}}}\int  f(k^{-1} x) |u|^{r(x)} \\
&=  \frac{1}{p_{+}} \int \left(  |\nabla u|^{p(x)}+|u|^{p(x)}\right)  - \frac{1}{{r_{-}}}\int   f(k^{-1} x) |u|^{r(x)}\\
& \quad- \frac{1}{q_{-}} \left( \int \left(  |\nabla u|^{p(x)}+|u|^{p(x)}\right) - \int  f(k^{-1} x) |u|^{r(x)} \right).
\end{align*}
Since $ p_{+} < q_{-} \leq q(x) \leq r(x) \ll p^{*}(x) $,
$$
J_{\lambda, k}(u) \geq \left( \frac{1}{p_{+}} - \frac{1}{q_{-}}\right) \int \left(  |\nabla u|^{p(x)}+|u|^{p(x)}\right),
$$
showing that $J$ is bounded from below and coercive on $ \mathcal{M}_{\lambda, k} $. \fim

\vspace{0.5 cm}

As an immediate consequence of the last lemma, we have

\begin{corollary}\label{ltdalem}
Let $ \{u_{n}\} $ be a sequence in $ \mathcal{M}_{\lambda, k} $ and  $ J_{\lambda, k}(u_{n}) \to c_{\lambda, k} $. Then $ \{u_{n}\} $ is bounded in $ W^{1,p(x)}(\mathbb{R}^{N}) $.
\end{corollary}

\vspace{0.5 cm}

The next lemma establishes that Nehari manifold $\mathcal{M}_{\lambda,k}$ has a positive distance from origin.

\begin{lem}\label{NehariDelta}
There exists  $\eta >0$  such that
\begin{equation} \label{R1}
\rho_1(u) \geq \eta,  \,\,\, \,\,\,\,\, \forall (u, \lambda, k)  \in \mathcal{M}_{\lambda, k} \times [0, \Lambda] \times \mathbb{N}.
\end{equation}
Moreover, if  $ E_{\lambda, k} (u) = J'_{\lambda, k} (u) u $,  we have that
\begin{equation} \label{E_lambda}
E'_{\lambda, k} (u) u \leq - \left( q_{-}-p_{+} \right) \eta
\,\,\,\,\, \forall (u, \lambda, k)  \in \mathcal{M}_{\lambda, k} \times [0, \Lambda] \times \mathbb{N}.
\end{equation}

\end{lem}
\noindent \textbf{Proof.}
Suppose by contradiction that (\ref{R1}) does not hold. Then, there is $ \{ u_{n} \} \subset \mathcal{M}_{\lambda, k} $ such that
$$
\rho_1(u_n) \to 0 ~~ \text{ as } n \to \infty,
$$
or equivalently, by Proposition \ref{modular2},
$$
\Vert u_{n} \Vert \to 0 \text{ as } n \to \infty.
$$
Since $ \{ u_{n} \} \subset \mathcal{M}_{\lambda, k} $ and $\|f\|_\infty \leq 1$, we derive
$$
\int  \left(  |\nabla u_{n}|^{p(x)}+|u_{n}|^{p(x)}\right)  \leq  \lambda \Vert g \Vert_{\infty} \int |u_{n}|^{q(x)} + \int |u_{n}|^{r(x)}.
$$
On the other hand, using the fact that $\Vert u_{n} \Vert < 1 $ for $n$ large enough, it follows from Propositions~{\ref{modular}} and \ref{modular2}
\begin{align*}
\Vert u_{n} \Vert^{p_{+}} \leq & \int   \left(  |\nabla u_{n}|^{p(x)}+|u_{n}|^{p(x)}\right)  \\
\leq & \;\lambda \Vert g \Vert_{\infty} \max\left\{ \Vert u_{n} \Vert^{q_{-}}_{q(x)} , \Vert u_{n} \Vert^{q_{+}}_{q(x)}  \right\} + \max\left\{ \Vert u_{n} \Vert^{r_{-}}_{r(x)} , \Vert u_{n} \Vert^{r_{+}}_{r(x)}  \right\}.
\end{align*}
By Sobolev embedding, there are positive constants $ c_{1} $ and $ c_{2} $ such that
\begin{align*}
\Vert u_{n} \Vert^{p_{+}} & \leq \lambda \Vert g \Vert_{\infty} c_{1} \max\left\{ \Vert u_{n} \Vert^{q_{-}} , \Vert u_{n} \Vert^{q_{+}}  \right\} +  c_{2} \max\left\{ \Vert u_{n} \Vert^{r_{-}} , \Vert u_{n} \Vert^{r_{+}}  \right\},
\end{align*}
and so, for $n$ large enough,
$$
\|u_n\|^{p_{+}} \leq \lambda c_1 \|g\|_{\infty} \|u_n\|^{q_{-}}+c_2 \|u_n\|^{r_{+}},
$$
obtaining an absurd, because $p_+ < q_{-} \leq r_-$. Therefore, (\ref{R1}) is proved.

\vspace{0.5 cm}

Next, we will show that (\ref{E_lambda}) occurs. For each $u \in \mathcal{M}_{\lambda, k}$, a simple calculus gives
$$
\begin{array}{l}
E'_{\lambda, k} (u) u \leq  \, p_{+} \displaystyle \int  \left(  |\nabla u|^{p(x)}+|u|^{p(x)}\right)   - \lambda q_{-} \int g(k^{-1}x)  |u|^{q(x)}
  - r_{-} \int f(k^{-1} x)   |u|^{r(x)} \\
  \mbox{} \\
\mbox{\hspace{1,4 cm} } \leq  \left( p_{+} - q_{-} \right)  \int  \left(  |\nabla u|^{p(x)}+|u|^{p(x)}\right)    + \left( q_{-} -r_{-} \right) \int f(k^{-1} x) |u|^{r(x)}.
\end{array}
$$
Since $ p_{+} < q_{-} \leq r_{-} $, it follows that
\[
E'_{\lambda, k} (u) u \leq -\left( q_{-}-p_{+} \right) \rho_1(u) \leq - \left( q_{-}- p_{+} \right) \eta, 
\]
finishing the proof. \fim

\vspace{0.5 cm}

As by product of the last lemma, we are able to prove that critical points of $J_{\lambda, k}$ restrict to  $\mathcal{M}_{\lambda, k}$ are in fact critical point of $J_{\lambda, k}$  on $ W^{1,p(x)}(\mathbb{R}^{N}) $.

\begin{lem}
If $ u_{0} \in \mathcal{M}_{\lambda, k} $ is a critical point of $ J_{\lambda, k} $ restricted to  $ {\mathcal{M}_{\lambda, k}} $, then $ u_{0} $ is a critical point of $ J_{\lambda, k}$ in $ W^{1,p(x)}(\mathbb{R}^{N}) $.
\end{lem}

\noindent \textbf{Proof.}
Once that $u_0$ is a critical point of  $ J_{\lambda, k} $ restricted to  $ {\mathcal{M}_{\lambda, k}} $, there is $ \tau \in \mathbb{R}^{N} $ such that
$$
J'_{\lambda, k}(u_{0}) = \tau E'_{\lambda, k} (u_{0}).
$$
By Lemma~\ref{NehariDelta}, we know that $ E'_{\lambda, k} (u_0) u_0 < 0 $, then we must have $\tau=0$. Thereby,
$$
J'_{\lambda, k}(u_{0}) = 0,
$$
implying that $ u_{0} $ is critical point of $ J_{\lambda, k}(u_{0})$ in  $ W^{1,p(x)}(\mathbb{R}^{N}) $. \fim

\vspace{0.5 cm}


The next result is very important in our arguments, because it implies that weak limit of $(PS)$ sequence is a critical point for the energy functional.

\begin{thm}\label{conv}
Let $ \{ u_{n} \} $ be a sequence in $  W^{1,p(x)}(\mathbb{R}^{N}) $ such that $ u_{n} \rightharpoonup u $ in $  W^{1,p(x)}(\mathbb{R}^{N}) $ and $ J'_{\lambda, k}(u_{n}) \to 0 $ as $ n \to \infty $. Then, for some subsequence, $ \nabla u_{n}(x) \to \nabla u(x) $ a.e. in $ \mathbb{R}^{N} $ as $ n \to \infty $ and $ J'_{\lambda, k}(u) = 0 $.
\end{thm}

\noindent \textbf{Proof.}
Let $ R > 0 $ and $ \phi \in C^{\infty}_{0}(\mathbb{R}^{N}) $ such that
$$ \phi = 0  \,\,\,\, \mbox{if} \,\,\,  |x| \geq 2R, \,\,\,  \,\,\,  \phi = 1 \,\,\, \mbox{if} \,\,\,  |x| \leq R \,\,\, \mbox{and} \,\,\, 0 \leq \phi(x) \leq 1 \,\,\,\forall x \in \mathbb{R}^{n}.
$$
In what follows, let us denote by $\{P_n\}$ the following sequence
$$
P_{n}(x) = \langle |\nabla u_{n}(x)|^{p(x)-2} \nabla u_{n}(x) - |\nabla u(x)|^{p(x)-2} \nabla u(x), \nabla u_{n}(x) - \nabla u (x) \rangle .
$$
From definition of $\{P_n\}$,
$$
\int_{B_{R}(0)} P_{n}  \leq  \int |\nabla u_{n}|^{p(x)} \phi - \int |\nabla u_{n}|^{p(x) - 2} \nabla u_{n} \nabla u \phi  - \int |\nabla u|^{p(x) - 2} \nabla u \nabla(u_{n} - u) \phi.
$$
Recalling that $u_n \rightharpoonup u$ in $W^{1,p(x)}(\mathbb{R}^{N})$, we have
\begin{align}
\int_{B_{R}(0)} |\nabla u|^{p(x) - 2} \nabla u \nabla(u_{n} - u)\phi \to 0 \quad \mbox{ as } n \to \infty,
\end{align}
and so,
\[
\int_{B_{R}(0)} P_{n} \leq \int |\nabla u_{n}|^{p(x)} \phi - \int |\nabla u_{n}|^{p(x) - 2} \nabla u_{n} \nabla u \phi + o_n(1).
\]
On the other hand, from $ J_{\lambda, k}'(u_{n})(\phi u_{n}) = o_n(1)$ and $ J_{\lambda, k}'(u_{n})(\phi u) = o_n(1)$,
\begin{align*}
\int_{B_{R}(0)} P_{n} & \leq o_n(1) -  \int |\nabla u_{n}|^{p(x) - 2}\nabla u_{n}\nabla \phi(u_{n} - u) \\
    & \quad - \int |u_{n}|^{p(x) - 2}u_{n}(u - u_{n}) \phi + \lambda \int g(k^{-1} x) |u_{n}|^{q(x) - 2}u_{n}(u_{n} -  u)\phi \\
        & \qquad + \int f(k^{-1} x) |u_{n}|^{r(x) - 2}u_{n}(u_{n} -  u)\phi.
\end{align*}
Thus,
\begin{align*}
\int_{B_{R}(0)} P_{n} & \leq o_n(1) +  c_{1}\int_{supt{\phi}} |\nabla u_{n}|^{p(x) - 1}|u_{n} - u| \\
    & \quad + c_1\int_{supt{\phi}} |u_{n}|^{p(x) - 1}|u_{n} - u| +c_1 \lambda \Vert g \Vert_{\infty} \int_{supt{\phi}} |u_{n}|^{q(x) - 1}|u_{n} -  u| \\
        & \qquad +c_1\int_{supt{\phi}} |u_{n}|^{r(x) - 1}|u_{n} -  u|.
\end{align*}
Combining H\"{o}lder's inequality and Sobolev embedding, we deduce that
\[
\int_{B_{R}(0)} P_{n} \to 0 \quad \mbox{as } n \to \infty.
\]
In what follows, let us consider the sets
\[
B_{R}^{+} = \left\{ x \in B_{R}(0) \, / \, p(x) \geq 2 \right\}\quad \text{ and } \quad B_{R}^{-} = \left\{ x \in B_{R}(0) \, / \, 1 < p(x) < 2 \right\}.
\]
Since
\[
P_{n}(x)  \geq \left\{
\begin{array}{ccc}
  \frac{2^{3- p_{+}}}{p_{+}} |\nabla u_{n} - \nabla u|^{p(x)} & \text{if} & p(x) \geq 2  \\ \\
  (p_{-} - 1)\frac{|\nabla u_{n} - \nabla u|^{2}}{\left(|\nabla u_{n}| + |\nabla u| \right)^{2 - p(x)}} &  \text{if} & 1 < p(x) < 2,
\end{array}
\right.
\]
we have
\begin{equation} \label{convE1}
\int_{B_{R}^{+}} |\nabla u_{n} - \nabla u|^{p(x)} dx \to 0 \text{ as} \,\, n \to \infty.
\end{equation}
Applying again H\"{o}lder's inequality,
$$
\hspace*{-2cm} \int_{B_{R}^{-}} |\nabla u_{n} - \nabla u|^{p(x)} \leq C \Vert g_{n}\Vert_{L^{\frac{2}{p(x)}}(B_{R}^{-})} \Vert h_{n}\Vert_{L^{\frac{2}{2 - p(x)}}(B_{R}^{-})},
$$
where
$$
g_{n}(x) = \frac{ |\nabla u_{n}(x) - \nabla u(x)|^{p(x)}}{\left( |\nabla u_{n}(x)| + |\nabla u(x)|\right)^{\frac{p(x)(2 - p(x))}{2}}}
$$
and
$$
h_{n}(x) = \left( |\nabla u_{n}(x)| + |\nabla u(x)| \right)^{\frac{p(x)(2 - p(x))}{2}}.
$$
By a direct computation, $ \{ \Vert h_{n}\Vert_{L^{\frac{2}{2 - p(x)}}(B_{R}^{-})}\}$ is a bounded sequence and
$$
\int_{B_{R}^{-}}|g_{n}|^{\frac{2}{p(x)}}  \leq C \int_{B_{R}^{-}}P_{n} .
$$
Then,
\begin{equation} \label{convE2}
\int_{B_{R}^{-}} |\nabla u_{n} - \nabla u|^{p(x)} \to 0 \text{ when } n \to \infty.
\end{equation}
From \eqref{convE1} and \eqref{convE2}, $ \nabla u_{n} \to \nabla u $ a.e. in $ B_{R}(0) $. Once that $R$ is arbitrary, it follows that for some subsequence
\[
\nabla u_{n}(x) \to \nabla u(x) \mbox{ a.e. in } \mathbb{R}^{N}.
\]
This combined with Lemma~\ref{Brezis-Lieb-2} gives
\[
|\nabla u_{n}|^{p(x) - 2} \nabla u_{n} \rightharpoonup |\nabla u|^{p(x) - 2} \nabla u \mbox{ in } (L^{p'(x)}(\mathbb{R}^{N}))^{N}.
\]
Now, using the fact that $J'_{\lambda,k}(u_n)v=o_n(1)$ for all $v \in W^{1,p(x)}(\mathbb{R}^{N}) $ together with the last limit, we derive that $J'_{\lambda,k}(u)v= 0 $ for all $v \in W^{1,p(x)}(\mathbb{R}^{N}) $, finishing the proof. \fim

\subsection{A result of compactness}

The next theorem is a version of a result compactness on Nehari manifolds  due to Alves \cite{alves05} for variable exponents. It establishes that problem (\ref{Poo}) has a ground state solution.

\begin{thm}\label{TeoComp}
Suppose that (\ref{p1}) holds and let $ \{u_{n}\} \subset \mathcal{M}_{\infty} $ be a sequence with $ J_{\infty}(u_{n}) \to c_{\infty} $. Then,
\begin{description}
\item[I.] $ u_{n} \to u $ in $ W^{1,p(x)}(\mathbb{R}^{N}) $,

or

\item[II.] There is $ \{y_{n}\} \subset \mathbb{Z}^{N}$ with $|y_n| \to +\infty$ and $w \in W^{1,p(x)}(\mathbb{R}^{N}) $ such that $ w_{n} = u_{n}(\cdot + y_{n}) \to w $ in $ W^{1,p(x)}(\mathbb{R}^{N}) $ and $J_{\infty}(w) = c_{\infty}$.
\end{description}
\end{thm}

\noindent \textbf{Proof.} Similarly to Corollary~\ref{ltdalem}, we can assume that $\{u_n\}$ is a bounded sequence, and so, there is $ u \in W^{1,p(x)}(\mathbb{R}^{N}) $ and a subsequence of $ \{ u_{n} \}$, still denoted by itself, such that $u_n  \rightharpoonup u $  in $ W^{1,p(x)}(\mathbb{R}^{N})$. Applying the Ekeland's variational principle, there is a sequence $ \{w_{n}\} $ in $ \mathcal{M}_{\infty} $ satisfying
\[
w_{n} = u_{n} + o_{n}(1), \quad J_{\infty}(w_{n}) \to c_{\infty}
\]
and
\begin{align} \label{eq1}
J'_{\infty}(w_{n}) - \tau_{n} E'_{\infty}(w_{n}) = o_{n}(1),
\end{align}
where $ (\tau_{n}) \subset \mathbb{R} $ and $ E_{\infty}(w) = J'_{\infty}(w)  w $, for any $ w \in W^{1,p(x)}(\mathbb{R}^{N}) $.

Since $ \{u_{n}\} \subset \mathcal{M}_{\infty} $, (\ref{eq1}) leads to
\[
\tau_{n} E'_{\infty}(w_{n})  w_{n} = o_{n}(1).
\]
By the arguments of Lemma~\ref{NehariDelta}, there exists $\delta >0$ such that
$$
|E'_{\infty}(w_{n})w_{n}| > \delta \,\,\, \forall n \in \mathbb{N}.
$$
From this, $ \tau_{n} \to 0 $ as $ n \to \infty $ and we can claim that
\[
J_{\infty}(u_{n}) \to c_{\infty} \,\,\, \mbox{and} \,\,\, J'_{\infty}(u_{n}) \to 0.
\]

Next, we will study the following  possibilities: $ u \neq 0 $ or $ u = 0 $.

\vspace{0.5 cm}

\noindent \textbf{Case 1:} $ u \neq 0 $.

\vspace{0.5 cm}

Similarly to Theorem~\ref{conv},  it follows that the below limits are valid for some subsequence:
\begin{itemize}
    \item $ u_{n}(x) \to u(x)$ \,\, and \,\, $\nabla u_{n}(x) \to \nabla u(x) $ a.e. in $ \mathbb{R}^{N}, $
    \item $ \displaystyle \int |\nabla u_{n}(x)|^{p(x)-2}\nabla u_{n}(x)\nabla v \to \int |\nabla u(x)|^{p(x)-2}\nabla u(x)\nabla v$,
    \item $  \displaystyle \int |u_{n}|^{p(x)-2} u_{n} v \to \int |u|^{p(x)-2} u_{n} v$,
\end{itemize}
and
\begin{itemize}
    \item $ \displaystyle \int |u_{n}|^{r(x) -2} u_{n} v \to \int |u|^{r(x) -2} u v  $
\end{itemize}
for any $ v \in W^{1,p(x)}(\mathbb{R}^{N})$. Consequently, $ u $ is critical point of  $ J_{\infty} $.
By Fatou's Lemma , it is easy to check that
\begin{align*}
c_{\infty} \leq & J_{\infty}(u) = J_{\infty}(u) - \frac{1}{r_{-}} J'_{\infty}(u) u \\
= &\int \left( \frac{1}{p(x)} - \frac{1}{r_{-}} \right) \left( | \nabla u |^{p(x)} + |u|^{p(x)} \right) + \int \left( \frac{1}{r_{-}} - \frac{1}{r(x)} \right)  |u|^{r(x)} \\
\leq & \liminf_{n \to \infty}\left\{ \int \left( \frac{1}{p(x)} - \frac{1}{r_{-}} \right) \left( | \nabla u_{n} |^{p(x)} + |u_{n}|^{p(x)} \right) \right.\\
&\qquad \left. + \int \left( \frac{1}{r_{-}} - \frac{1}{r(x)} \right)  |u_{n}|^{r(x)}\right\} \\
= & \liminf_{n \to \infty} \left\{ J_{\infty}(u_{n}) - \frac{1}{r_{-}} J'_{\infty}(u_{n}) u_{n} \right\} = \, c_{\infty} .\\
\end{align*}
Hence,
\begin{align*}
\lim_{n \to \infty} \int \left( | \nabla u_{n} |^{p(x)} + |u_{n}|^{p(x)} \right) = \int \left( | \nabla u |^{p(x)} + |u|^{p(x)} \right),
\end{align*}
implying that $ u_{n} \to u $ in $ W^{1,p(x)}(\mathbb{R}^{N})$.

\vspace{0.5 cm}

\noindent \textbf{Case 2:} $ u = 0 $.

\vspace{0.5 cm}

In this case, we claim that there are $R, \xi>0$ and $ \{ y_{n} \} \subset \mathbb{R}^{N} $ satisfying
\begin{align}\label{lionsfalse}
\limsup_{n \to \infty} \int_{B_{R}(y_{n})} |u_{n}|^{p(x)} \geq \xi.
\end{align}
If the claim is false, we must have
\begin{align*}
\limsup_{n \to \infty} \sup_{y \in \mathbb{R}^{N}} \int_{B_{R}(y)} |u_{n}|^{p(x)} = 0.
\end{align*}
Thus, by a Lions-type result for variable exponent proved in \cite[Lema~3.1]{Fan2001},
$$
u_{n} \to 0  \mbox{ in }  L^{s(x)}(\mathbb{R}^{N}),
$$
for any $ s \in C(\mathbb{R}^{N}) $ with $ p \ll s \ll p^{*}$.

 Recalling $ J'_{\infty}(u_{n}) u_{n} = o_{n}(1) $, the last limits yield
\[
\int \left( | \nabla u_{n} |^{p(x)} + |u_{n}|^{p(x)} \right) = o_{n}(1),
\]
or equivalently
$$
u_n \to 0 \,\,\, \mbox{in} \,\,\, W^{1,p(x)}(\mathbb{R}^{N}),
$$
leading to $c_{\infty} = 0 $, which is absurd. This way, (\ref{lionsfalse}) is true. By a routine argument, we can assume that $ {y}_{n} \in \mathbb{Z}^{N} $ and $ |y_{n}| \to \infty $ as $ n \to \infty $.
Setting
$$
w_{n}(x) = u_{n}(x + {y}_{n}),
$$
and using the fact that $ p $ and $ r $ are $ \mathbb{Z}^{N}$-periodic, a change of variable gives
$$
J_{\infty}(w_{n}) = J_{\infty}(u_{n}) \,\,\, \mbox{and} \,\,\, \|J'_{\infty}(w_{n})\| = \|J'_{\infty}(u_{n})\|,
$$
showing that $ \{ w_{n} \} $ is a sequence $(PS)_{c_{\infty}} $ for $ J_{\infty} $.  If $ w \in  W^{1,p(x)}(\mathbb{R}^{N}) $ denotes the weak limit of $ \{ w_{n} \} $, it follows from (\ref{lionsfalse}),
$$
\int_{B_{{R}}(0)} |w|^{p(x)} \geq \xi,
$$
showing that $ w \neq 0 $.

Repeating the same argument of the first case for the sequence $ \{w_{n}\} $, we deduce that $ w_{n} \to w $ in $ W^{1,p(x)}(\mathbb{R}^{N}) $, $ w \in \mathcal{M}_{\infty} $  and $ J_{\infty}(w) = c_{\infty} $.  \fim

\subsection{Estimates involving the minimax levels}

The main goal of this section is to prove some estimates involving the minimax levels $c_{\lambda, k},c_{0, k}$ and
$c_{\infty}$.

First of all, we recall the inequalities
$$
J_{\lambda, k} (u) \leq J_{0, k}(u) \,\,\, \mbox{and} \,\,\, J_{\infty}(u) \leq J_{0, k}(u) \,\,\,\,\, \forall u \in W^{1,p(x)}(\mathbb{R}^{N}),
$$
which imply
\[
c_{\lambda, k} \leq c_{0, k}\quad  \mbox{ and }  \quad c_{\infty} \leq c_{0, k}.
\]

\begin{lem}\label{c0<cf00}
The minimax levels $c_{0, k}$ and $c_{f_{\infty}}$ satisfy the inequality \linebreak $c_{0, k} < c_{f_{\infty}}$. Hence, $c_{\infty} < c_{f_{\infty}}$.
\end{lem}
\noindent \textbf{Proof.}
In a manner analogous to Theorem~\ref{TeoComp}, there is $ U \in W^{1,p(x)}(\mathbb{R}^{N}) $ verifying
\[
J_{f_{\infty}}(U) = c_{f_{\infty}} \quad \mbox{ and } \quad J'_{f_{\infty}}(U) = 0.
\]
From Lemma~3.6 in \cite{Fan2008}, there exists $ t > 0 $ such that $ t U \in \mathcal{M}_{0, k} $. Thus,
$$
c_{0, k} \leq  J_{0, k}(tU) = \int \frac{t^{p(x)}}{p(x)} \left( |\nabla U|^{p(x)} + |U|^{p(x)}\right) - \int  f(k^{-1} x) \frac{t^{r(x)}}{r(x)}  |U|^{r(x)}.
$$
Since that by $(H2)$, $f_{\infty} < f(x)$ for all $x \in \mathbb{R}^{N}$, we derive
$$
c_{0, k} <  J_{f_{\infty}}(tU) \leq \max_{s \geq 0}J_{f_{\infty}}(sU) = J_{f_{\infty}}(U) = c_{f_{\infty}}.
$$
\fim

\vspace{0.5 cm}

Using the last lemma, we are able to prove that  $ J_{\lambda, k} $  verifies the $(PS)_{d}$ condition for some values of $d$.

\begin{lem}\label{Cond-PS}
The functional $ J_{\lambda, k} $ satisfies the $(PS)_{d}$ condition for \linebreak $ d \leq c_{\infty} + \varrho $, where $\varrho=\frac{1}{2}(c_{f_\infty}-c_\infty)>0$.
\end{lem}
\begin{pf}
Let $ \{v_{n}\} \subset W^{1,p(x)}(\mathbb{R}^{N}) $ be a $(PS)_{d}$ sequence for functional $ J_{\lambda, k} $ with $ d \leq c_{\infty} + \varrho$. Similarly to Corollary~\ref{ltdalem}, $ \{v_{n}\} $ is a bounded sequence in $ W^{1,p(x)}(\mathbb{R}^{N}) $, and so, for some subsequence, still denoted by $ \{v_{n}\} $,
\[
v_{n} \rightharpoonup v \mbox{ in } W^{1,p(x)}(\mathbb{R}^{N}),
\]
for some $v \in W^{1,p(x)}(\mathbb{R}^{N}).$ Now, we claim that
\begin{align}
J_{\lambda, k} (v_{n}) - J_{0, k} (w_{n}) - J_{\lambda, k}(v) = o_{n}(1) \label{J-J0-on1}
\end{align}
and
\begin{align}
\Vert J'_{\lambda, k} (v_{n}) - J'_{0, k} (w_{n}) - J'_{\lambda, k}(v) \Vert = o_{n}(1), \label{J'-J0'-on1}
\end{align}
where $ w_{n} = v_{n} - v $.

Indeed, proceeding as in proof of Theorem ~\ref{conv}, we have the following convergences
\[
\nabla v_{n}(x) \to \nabla v(x) \,\,\, \mbox{and} \,\,\, v_{n}(x) \to v(x) \,\,\, \mbox{ a.e. in} \,\,\, \mathbb{R}^{N}.
\]

Applying Lemma ~\ref{Brezis-Lieb-1}, it follows that
\begin{align*}
J_{\lambda, k}(v_{n})  = J_{0,k}(w_{n}) + J_{\lambda, k}(v) + o_{n}(1), \label{eq7}
\end{align*}
showing~(\ref{J-J0-on1}). The equality (\ref{J'-J0'-on1}) follows combining \ref{H2} with Lemmas \ref{Brezis-Lieb-2} and \ref{Brezis-Lieb-3}.

Since $ J'_{\lambda, k}(v) = 0 $ and  $ J_{\lambda, k}(v) \geq 0 $, from (\ref{J-J0-on1})-(\ref{J'-J0'-on1}), we have that $w_{n} = v_{n} - v$ is a $(PS)_{d^{*}}$ sequence for $J_{0, k}$ with  $d^*=d - J_{\lambda, k}(v) \leq  c_{\infty} +\varrho$.

\bigskip

\begin{claim} \label{C2} There is $ R > 0 $ such that
\[
\limsup_{n \to \infty} \sup_{y \in \mathbb{R}^{N}} \int_{B_{R}(y)}|w_{n}|^{p(x)} = 0.
\]
\end{claim}
If the claim is true, we have
\[
\int |w_{n}|^{r(x)} \to 0.
\]
On the other hand, by (\ref{J'-J0'-on1}),  we know that $ J'_{0, k}(w_{n}) = o_{n}(1) $, then

\begin{align*}
\int \left( | \nabla w_{n}|^{p(x)} + |w_{n}|^{p(x)}\right) = o_{n}(1),
\end{align*}
showing that  $ w_{n} \to 0 $ in $ W^{1,p(x)}(\mathbb{R}^{N}) $, and so, $ v_{n} \to v $ in $ W^{1,p(x)}(\mathbb{R}^{N}). $

\bigskip

\noindent\textbf{Proof of Claim \ref{C2}:}
If the claim is not true, for each $ R > 0 $ given, we find $ \xi > 0 $ and $ \{y_{n}\} \subset \mathbb{Z}^{N} $  verifying
\[
\limsup_{n \to \infty}\int_{B_{R}(y_{n})} |w_{n}|^{p(x)} \geq \xi > 0.
\]
Once that $w_n \rightharpoonup 0$ in $W^{1,p(x)}(\mathbb{R}^{N}) $, it follows that $\{y_n\}$ is an unbounded sequence.  Setting
\[
\tilde{w}_{n} = w_{n}(\cdot + y_{n}),
\]
we have that $ \{\tilde{w}_{n}\} $ is also a $(PS)_{d^*}$ sequence for $J_{0, k}$, and so, it must be bounded. Then, there are $ \tilde{w} \in W^{1,p(x)}(\mathbb{R}^{N})  $ and a subsequence of $ \{\tilde{w}_{n}\} $, still denoted by itself,  such that
$$
\tilde{w}_{n} \rightharpoonup \tilde{w} \in W^{1,p(x)}(\mathbb{R}^{N})\setminus\{0\}.
$$
Moreover, since $ J'_{0, k}(w_{n}) \phi( \cdot - y_{n}) = o_{n}(1) $ for each $ \phi \in  W^{1,p(x)}(\mathbb{R}^{N})$ and $ \nabla \tilde{w}_{n}(x) \to \nabla \tilde{w}(x) $ a.e. in $ \mathbb{R}^{N} $, we obtain
\begin{align*}
\int \left( |\nabla \tilde{w}|^{p(x) - 2} \nabla \tilde{w} \nabla \phi + | \tilde{w} |^{p(x) - 2} \tilde{w} \phi \right) = \int f_{\infty} | \tilde{w} |^{r(x) - 2} \tilde{w} \phi,
\end{align*}
from where it follows that $ \tilde{w} $ is a weak solution of the Problem~(\ref{Pf00}). Consequently, after some routine calculations, we get
$$
c_{f_{\infty}}  \leq J_{f_{\infty}}(\tilde{w}) = J_{f_{\infty}}(\tilde{w}) - \frac{1}{r_{-}} J'_{f_{\infty}}(\tilde{w})\tilde{w} \leq \liminf_{n \to \infty} \left\{J_{0, k}(w_{n}) - \frac{1}{r_{-}} J'_{0, k}(w_{n})  w_{n}\right\} = d^*
$$
implying that $c_{f_{\infty}}  \leq c_{\infty} + \varrho$, which is an absurd because $\varrho < c_{f_{\infty}}  - c_{\infty}$. Therefore, the Claim \ref{C2} is true.
\end{pf}

In what follows, let us fix $ \rho_{0}, r_{0} > 0 $ satisfying
\begin{itemize}
  \item $ \overline{B_{\rho_{0}}(a_{i})} \cap \overline{B_{\rho_{0}}(a_{j})} = \emptyset $ for $ i \neq j$ \,\,\, \mbox{and} \,\,\, $i,j \in \{1,...,\ell\}$
  \item $ \bigcup^{\ell}_{i = 1}B_{\rho_{0}}(a_{i}) \subset B_{r_0}(0) $.

  \item $K_{\frac{\rho_{0}}{2}} = \bigcup^{\ell}_{i = 1}\overline{B_{\frac{\rho_{0}}{2}}(a_{i})}$
\end{itemize}
Besides this, we define the function $ Q_{k} : W^{1, p(x)}(\mathbb{R}^{N}) \to \mathbb{R} $ by
\begin{align*}
Q_{k}(u) = \frac{\int \chi (k^{-1} x)|u|^{p_{+}}}{\int |u|^{p_{+}}},
\end{align*}
where  $ \chi : \mathbb{R}^{N}  \to \mathbb{R}^{N} $ is given by
\[
\chi(x) = \left\{
\begin{array}{ccc}
  x & \mbox{if} & |x| \leq r_{0} \\
  r_{0} \frac{x}{|x|} & \mbox{if} & |x| > r_{0}.
\end{array}
\right.
\]

The next two lemmas will be useful to get important $(PS)$-sequences associated with $ J_{\lambda, k} $.

\begin{lem}\label{lemK}
There are $ \delta_{0} > 0 $ and $ k_1 \in \mathbb{N} $ such that if $ u \in \mathcal{M}_{0, k} $ and $ J_{0, k}(u) \leq c_{\infty} + \delta_{0} $, then
\[
Q_{k}(u) \in K_{\frac{\rho_{0}}{2}} \,\,\,\, \mbox{for} \,\,\, k \geq  k_1.
\]
\end{lem}
\noindent \textbf{Proof.}
If the lemma does not occur, there must be $ \delta_{n} \to 0 $, $ k_{n} \to +\infty $ and $ u_{n} \in \mathcal{M}_{0, k_{n}} $ satisfying
\[
J_{0, k_{n}}(u_{n}) \leq c_{\infty} + \delta_{n}
\]
and
\[
Q_{k_{n}}(u_{n}) \not\in K_{\frac{\rho_{0}}{2}}.
\]
Fixing $ s_{n} > 0 $ such that $ s_{n} u_{n} \in \mathcal{M}_{\infty} $, we have that
\[
c_{\infty} \leq J_{\infty}(s_{n} u_{n}) \leq J_{0, k_{n}} (s_{n}u_{n}) \leq \max_{t \geq 0 } J_{0, k_{n}} (tu_{n}) = J_{0, k_{n}}(u_{n}) \leq c_{\infty} + \delta_{n}
\]
hence,
\[\{s_{n} u_{n}\} \subset \mathcal{M}_{\infty} \,\,\,\, \mbox{and} \,\,\,\, J_{\infty}(s_{n} u_{n}) \to c_{\infty}.
\]

Applying the variational principle of Ekeland, we can assume without loss of generality that
$\{s_{n} u_{n}\} \subset \mathcal{M}_{\infty} $  is a sequence $ (PS)_{c_{\infty}} $ for $ J_{\infty} $, that is,
$$
J_\infty(s_n u_n) \to c_\infty \,\,\,\, \mbox{and} \,\,\,\, J'_{\infty}(s_n u_n) \to 0.
$$
According to Theorem \ref{TeoComp}, we must consider the ensuing cases:
\begin{description}
  \item[i)] $ s_{n}u_{n} \to U \neq 0 $ in $ W^{1,,p(x)}(\mathbb{R}^{N}) $; \par
\end{description}
or
\begin{description}
\item[ii)] There exists $ \{y_{n}\} \subset \mathbb{Z}^{N} $ with $|y_n| \to +\infty $ such that $ v_{n} = s_{n}u(\cdot + y_{n}) $ is convergent in $ W^{1,,p(x)}(\mathbb{R}^{N}) $ for some $ V \in  W^{1,p(x)}(\mathbb{R}^{N}) \setminus \{0\}$.
\end{description}

By a direct computation, we can suppose that $ s_{n} \to s_{0} $ for some $ s_{0} > 0 $. Therefore, without loss of generality,
we can assume that
$$
u_{n} \to U  \,\,\, \mbox{or} \,\,\,\, v_{n} = u( \,\, \cdot + y_{n}) \to  V \,\,\,\, \mbox{in} \,\,\,  W^{1, p(x)}(\mathbb{R}^{N}).
$$

\bigskip

\noindent\textbf{Analysis of} $\mathbf{i)}$.

\bigskip

By Lebesgue's dominated convergence theorem
\[
Q_{k_{n}}(u_{n}) = \frac{\int \chi({k_{n}}^{-1} x)|u_{n}|^{p_{+}}}{\int |u_{n}|^{p_{+}}} \to \frac{\int \chi(0)|U|^{p_{+}}}{\int |U|^{p_{+}}} = 0 \in K_{\frac{\rho_{0}}{2}},
\]
implying $ Q_{k_{n}}(u_{n}) \in K_{\frac{\rho_{0}}{2}} $ for $n $ large, which is an absurd.

\bigskip

\noindent\textbf{Analysis of} $\mathbf{ii)}$.

\bigskip
Using again the Ekeland's variational principle, we can suppose that $ J'_{0, k_{n}}(u_{n}) = o_{n}(1) $. Hence, $ J'_{0, k_{n}} (u_{n}) \phi(\cdot - y_{n}) = o_{n}(1) $ for  any $ \phi \in W^{1,,p(x)}(\mathbb{R}^{N})$, and so,
\begin{equation}
o_{n}(1)  = \int \left(  |\nabla v_{n}|^{p(x) - 2} \nabla v_{n} \nabla \phi + |v_{n}|^{p(x) - 2} v_{n} \phi\right) - \int f({k_{n}}^{-1} (x + y_{n}))  |v_{n}|^{r(x) - 2}v_{n} \phi. \label{eq5}
\end{equation}
The last limit implies that for some subsequence,
$$
\nabla v_{n}(x) \to \nabla V(x) \,\,\, \mbox{and} \,\,\,  v_{n}(x) \to  V(x) \,\,\, \mbox{a.e in} \,\,\, \mathbb{R}^{N}.
$$
Now, we will study two cases:
\begin{description}
  \item[I)] $ |{k_{n}}^{-1}y_{n}| \to +\infty $
\end{description}
and
\begin{description}
  \item[II)] $ {k_{n}}^{-1}y_{n} \to y $, for some $y \in \mathbb{R}^{N}$.
\end{description}

If I) holds, it follows that
\[
\int \left(  |\nabla V|^{p(x) - 2} \nabla V \nabla \phi + |V|^{p(x) - 2} V \phi\right) = \int f_{\infty} |V|^{r(x) - 2}  V \phi,
\]
showing that $ V $ is a nontrivial weak solution of the problem~(\ref{Pf00}). Now, by Fatou's Lemma,
$$
c_{f_{\infty}} \leq J_{f_{\infty}}(V)  = J_{f_{\infty}}(V) - \frac{1}{r_{-}}J'_{f_{\infty}}(V)V \leq \liminf_{n \to \infty} \left\{  J_{\infty}(u_{n}) - \frac{1}{r_{-}}J'_{\infty}(u_{n})u_{n}  \right\} = c_{\infty},
$$
or equivalently, $ c_{f_{\infty}} \leq c_{\infty} $, contradicting the Lemma~\ref{c0<cf00}.

\bigskip



Now, if  $ {k_{n}}^{-1}y_{n} \to y $ for some $y \in \mathbb{R}^{N}$, then $ V $ is a weak solution of the following problem
\begin{align}
\left\{
\begin{array}
[c]{rcl}%
-\Delta_{p(x)} u + |u|^{p(x) - 2} u & = & f(y)|u|^{r(x) - 2}u  , \quad
\mathbb{R}^{N}\\
\mbox{}\\
u \in W^{1, p(x)}(\mathbb{R}^{N}). &  &
\end{array}
\right. \tag{$ P_{f(y)} $}\label{Pfy}%
\end{align}
Repeating the previous argument, we deduce that
\begin{align}
c_{f(y)} \leq c_{\infty}, \label{eq6}
\end{align}
where $ c_{f(y)} $ the mountain pass level of the functional  $ J_{f(y)} : W^{1,p(x)}(\mathbb{R}^{N}) \to \mathbb{R} $ given by
\begin{align*}
J_{f(y)}(u)   =\int \frac{1}{p(x)}\left(  |\nabla
u|^{p(x)}+|u|^{p(x)}\right)    -  \int \frac{ f(y) }{r(x)} |u|^{r(x)}.
\end{align*}
Observe that
$$
c_{f(y)} = \inf_{u \in \mathcal{M}_{f(y)}} J_{f(y)}(u)
$$
where
$$
\mathcal{M}_{f(y)} = \left\{ u \in W^{1,p(x)}(\mathbb{R}^{N})\setminus\{ 0\}: J'_{f(y)}(u) u = 0 \right\}.
$$
If $ f(y) < 1 $, a similar argument explored in the proof of Lemma~\ref{c0<cf00} shows that $ c_{f(y)} > c_{\infty} $, contradicting the inequality (\ref{eq6}). Thereby, $ f(y) = 1$ and $ y = a_{i} $ for some $ i = 1, \cdots \ell $.
Hence,
\begin{align*}
Q_{k_{n}}(u_{n})  &=  \frac{\int \chi({k_{n}}^{-1} x)|u_{n}|^{p_{+}}}{\int |u_{n}|^{p_{+}}}\\
 & =  \frac{\int \chi({k_{n}}^{-1} x + {k_{n}}^{-1} y_{n})|v_{n}|^{p_{+}}}{\int |v_{n}|^{p_{+}}} \to  \frac{\int \chi(y)|V|^{p_{+}}}{\int |V|^{p_{+}}}=a_i \in K_{\frac{\rho_{0}}{2}}, \\
\end{align*}
implying that $ Q_{k_{n}}(u_{n}) \in K_{\frac{\rho_{0}}{2}} $ for $ n $ large, which is a contradiction, since by assumption $ Q_{k_{n}}(u_{n}) \not\in K_{\frac{\rho_{0}}{2}} $.
\fim

\begin{lem}\label{lemK2}
Let $ \delta_{0}>0 $ given in Lemma \ref{lemK} and $ k_3=\max\{k_1,k_2\}$. Then, there is $ \Lambda^{*} > 0 $ such that
\[
Q_{k}(u) \in K_{\frac{\rho_{0}}{2}}, \,\,\,\,\,\, \forall (u, \lambda, k) \in \mathcal{A}_{\lambda, k} \times [0, \Lambda_*) \times ( [k_3,+\infty) \cap \mathbb{N}).
\]
\end{lem}

\noindent \textbf{Proof.}
Observe that
$$
J_{\lambda, k}(u) = J_{0, k}(u) - \lambda \int \frac{g({k}^{-1} x)}{q(x)} |u|^{q(x)} \,\,\, \forall u \in W^{1,p(x)}(\mathbb{R}^{N}).
$$
In what follows, let $ t_{u} > 0 $ such that $ t_{u}u \in \mathcal{M}_{0, k} $. Then,
\begin{align}
J_{0, k}(t_{u} u) & = J_{\lambda, k}(t_{u}u) + \lambda \int \frac{g({k}^{-1} x)}{q(x)} (t_{u})^{q(x)} |u|^{q(x)} \nonumber \\
    & \leq  \max_{t \geq 0} J_{\lambda, k}(tu)+ \lambda \int \frac{g({k}^{-1} x)}{q(x)}  (t_{u})^{q(x)} |u|^{q(x)}  \label{Z1}.
\end{align}

\begin{claim}\label{tuNerahi} \mbox{}\\

\noindent {\bf a)} There is a constant $ R > 0 $ such that
\[
\mathcal{A}_{\lambda, k} = \left\{ u \in \mathcal{M}_{\lambda, k}; J_{\lambda, k}(u) < c_{\infty} + \frac{\delta_{0}}{2} \right\} \subset B_{R}(0),
\]
for $k \geq k_1$,  that is, $ \mathcal{A}_{\lambda, k} $ is bounded set, where $k_1$ was given in Lemma \ref{lemK}. Moreover, $ R $ is independent of $ \lambda $ and $k$.

\vspace{0.5 cm}

\noindent {\bf b)} Let $ u \in \mathcal{A}_{\lambda, k} $ and $ t_{u} > 0 $ such that $ t_{u} u \in \mathcal{M}_{0, k} $. Then, given $\Lambda >0$, there are $ C > 0 $ and $k_2 \in \mathbb{N}$ such that
\[
0 \leq t_{u} \leq C, \quad \mbox{ for all } (u, \lambda, k)  \in \mathcal{A}_{\lambda, k} \times [0,\Lambda] \times ([ k_2, +\infty) \cap \mathbb{N}).
\]
\end{claim}

\noindent {\bf Proof of a):} Let $ u \in \mathcal{M}_{\lambda, k} $ such that $ J_{\lambda, k}(u) < c_{\infty} + \frac{\delta_{0}}{2}$ for $ k \geq k_1 $. Then,
\begin{align*}
\int \left( |\nabla u|^{p(x)} + |u|^{p(x)} \right) - \lambda \int g({k}^{-1} x) |u|^{q(x)} - \int f({k}^{-1} x)|u|^{r(x)} = 0
\end{align*}
and
\begin{align*}
\int  \frac{1}{p(x)} \left( |\nabla u|^{p(x)} + |u|^{p(x)} \right) - & \lambda \int  \frac{g({k}^{-1} x)}{q(x)} |u|^{q(x)} - \int   \frac{f({k}^{-1} x)}{r(x)} |u|^{r(x)}  < c_{\infty} + \frac{\delta_{0}}{2}.
\end{align*}
Combining the last two expressions, we obtain
\begin{align*}
\left( \frac{1}{p_{+}} - \frac{1}{q_{-}} \right) \int  \left( |\nabla u|^{p(x)} + |u|^{p(x)} \right) + \left( \frac{1}{q_{-}} - \frac{1}{r_{-}} \right) \int f({k}^{-1} x) |u|^{r(x)}
 < c_{\infty} + \frac{\delta_{0}}{2}.
\end{align*}
Therefrom,
\begin{align*}
\int \left( |\nabla u|^{p(x)} + |u|^{p(x)} \right)  < (c_{\infty} + \frac{\delta_{0}}{2})\left( \frac{1}{p_{+}} - \frac{1}{q_{-}} \right)^{-1},
\end{align*}
proving a).

\vspace{0.5 cm}

\noindent {\bf Proof of b):} Supposing by contradiction that the lemma does not hold, there is $ \{u_{n}\} \subset \mathcal{A}_{\lambda_n, k_n} $ with $\lambda_n \to 0$ and $k_n \to +\infty$ such that $ t_{u_{n}} u_{n} \in \mathcal{M}_{0, k_n} $ and $ t_{u_{n}} \to \infty $ as $ n \to \infty $. Without loss of generality, we assume that $ t_{u_{n}} \geq 1$. As $ t_{u_{n}} u_{n} \in \mathcal{M}_{0, k_n} $, we derive
$$
(t_{u_{n}})^{p_+} \int  \left( |\nabla u_{n}|^{p(x)} + |u_{n}|^{p(x)} \right) \geq f_{\infty} (t_{u_{n}})^{r_{-}}  \int |u_{n}|^{r(x)},
$$
or equivalently,
\begin{align}
\int \left( |\nabla u_{n}|^{p(x)} + |u_{n}|^{p(x)} \right) \geq  f_{\infty} t_{u_{n}}^{r_{-}- p_{+}}  \int |u_{n}|^{r(x)} .  \label{modular-tu}
\end{align}
Now, we claim that there is $ \eta_{1} > 0 $ such that
\begin{equation} \label{eta1}
\int |u_n|^{r(x)} > \eta_{1} \,\,\, \forall n \in \mathbb{N}.
\end{equation}
Indeed, arguing by contradiction, if $ \int  |u_{n}|^{r(x)} \to 0 $, by interpolation it follows that
$ \int_{\mathbb{R}^{N}} |u_n|^{q(x)} \to 0$ . Since $ u_{n} \in \mathcal{M}_{\lambda_{n}, k_{n}} $,
$$
\int \left( |\nabla u_{n}|^{p(x)} + |u_{n}|^{p(x)} \right) \leq \lambda_n \Vert g \Vert_{\infty} \int |u_{n}|^{q(x)} +  \int |u_{n}|^{r(x)} = o_{n}(1),
$$
or equivalently,
$$
u_{n}  \to 0 \,\,\, \mbox{in} \,\,\, W^{1, p(x)}(\mathbb{R}^{N}),
$$
which contradicts Lemma~\ref{NehariDelta}, proving (\ref{eta1}). Thereby, from inequality (\ref{modular-tu}),
$$
\rho_{1}(u_{n})=\int \left( |\nabla u_{n}|^{p(x)} + |u_{n}|^{p(x)} \right)  \to +\infty,
$$
implying that $ \{u_{n}\}$ is a unbounded sequence. However, this is impossible, because by item a), $ \{u_{n} \} $ is bounded, showing that b) holds.

Now, combining Claim \ref{tuNerahi}-b) with (\ref{Z1}), we get
\begin{align*}
J_{0, k}(t_{u} u) \leq J_{\lambda, k}(u) + \frac{\lambda}{q_{-}}\Vert g \Vert_{\infty} C^{q_+} \int  |u|^{q(x)}.
\end{align*}
Once that $u \in \mathcal{A}_{\lambda, k}$, we derive
\begin{align*}
J_{0, k}(t_{u} u) < c_{\infty} + \frac{\delta_{0}}{2} + \lambda c_{2}  \int|u|^{q(x)}.
\end{align*}
Using the Sobolev embedding combined with Claim \ref{tuNerahi}-a), we obtain
\[
J_{0, k} (t_{u} u) < c_{\infty} + \frac{\delta_{0}}{2} + c_{3} \lambda \quad \forall u \in \mathcal{A}_{\lambda, k}
\]
where $c_{3} $ is a positive constant. Setting $ \Lambda^{*} : = {\delta_{0}}/{2c_{3}}$  and $\lambda \in [0, \Lambda^*)$, it follows that
\[
t_{u}u \in \mathcal{M}_{0, k} \quad \mbox{ and } \quad J_{0, k} (t_{u} u) < c_{\infty} + \delta_{0}.
\]
Then, by Lemma~\ref{lemK},
\[
Q_{k}(t_{u} u) \in K_{\frac{\rho_{0}}{2}}.
\]
Now, it remains to note that
\[
Q_{k}(u) = Q_{k}(t_{u} u),
\]
to conclude the proof of lemma. \fim

\vspace{0.5 cm}

From now on, we will use the ensuing notation
\begin{itemize}
  \item $ \theta^{i}_{\lambda, k} = \left\{u \in \mathcal{M}_{\lambda, k} ; |Q_{k}(u) - a_{i} | < \rho_{0} \right\}$,
  \item $ \partial\theta^{i}_{\lambda, k} = \left\{u \in \mathcal{M}_{\lambda, k} ; |Q_{k}(u) - a_{i} | = \rho_{0} \right\}$,
  \item $ \beta^{i}_{\lambda, k} = \displaystyle\inf_{u \in \theta^{i}_{\lambda, k}} J_{\lambda, k}(u) $
  \end{itemize}
and
\begin{itemize}
 \item  $ \tilde{\beta}^{i}_{\lambda, k} = \displaystyle\inf_{u \in \partial\theta^{i}_{\lambda, k}} J_{\lambda, k}(u) .$
\end{itemize}

The above numbers are very important in our approach, because we will prove that there is a $(PS)$ sequence of $J_{\lambda, k}$ associated with each $\theta^{i}_{\lambda, k}$ for $i=1,2,...,\ell$. To this end, we need of the following technical result

\begin{lem} \label{rho} There is $ k^{*} \in \mathbb{N} $ such that
$$
\beta^{i}_{\lambda, k}  < c_{\infty} + \varrho \,\,\, \mbox{and} \,\,\, \beta^{i}_{\lambda, k} < \tilde{\beta}^{i}_{\lambda, k},
$$
for all $ \lambda \in [0, \Lambda^*)$ and $ k \geq k^{*}$, where $\varrho=\frac{1}{2}(c_{f_\infty}-c_\infty)>0.$
\end{lem}
\noindent \textbf{Proof.} From now on,  $ U \in W^{1, p(x)}(\mathbb{R}^{N}) $ is a ground state solution associated with (\ref{Poo}), that is,
\[
J_{\infty}(U) = c_{\infty} \quad \mbox{ and } \quad  J'_{\infty}(U) = 0 \,\,\, \mbox{( See Theorem \ref{TeoComp} )}.
\]
For $ 1 \leq i \leq \ell$ and $ k \in \mathbb{N} $, we define the function $ \widehat{U}^{i}_{k} : \mathbb{R}^{N} \to \mathbb{R}$ by
\[
\widehat{U}^{i}_{k}(x) = U( x - ka_{i}).
\]

\begin{claim}\label{ltda-Jl-coo} For all $i \in \{1,...,\ell \}$, we have that
\[
\limsup_{k \to +\infty}(\sup_{t \geq 0 } J_{\lambda, k}(t\widehat{U}^{i}_{k})) \leq c_{\infty}.
\]
\end{claim}
Indeed, since $p,q$ and $ r $ are  $ \mathbb{Z}^{N}$-periodic, and $ a_{i} \in \mathbb{Z}^{N}$, a change variable gives
\begin{align*}
J_{\lambda, k}(t\widehat{U}^{i}_{k}) = &\int  \frac{t^{p(x)}}{p(x)} \left( |\nabla U|^{p(x)} + |U|^{p(x)} \right)
- \lambda  \int g(k^{-1} x + a_{i})\frac{t^{q(x)}}{q(x)} \left| U \right|^{q(x)}\\
& \qquad - \int f(k^{-1} x + a_{i}) \frac{t^{r(x)}}{r(x)} |U|^{r(x)}.
\end{align*}
Moreover, we know that there exists $ s = s(k)  > 0 $ such that
\begin{align*}
\max_{t \geq 0} J_{\lambda, k}(t\widehat{U}^{i}_{k}) = J_{\lambda, k}(s\widehat{U}^{i}_{k}).
\end{align*}
By a direct computation, it follows that $ s(k) \not\to 0 $ and  $ s(k) \not\to \infty $ as $ k \to \infty $. Thus, without loss of generality, we can assume $ s(k) \to s_{0}>0$ as $ k \to \infty $. Thereby,
\begin{align*}
\lim_{k \to \infty} \left( \max_{t \geq 0}J_{\lambda, k}(t\widehat{U}^{i}_{k})  \right) &= \int \frac{s_{0}^{p(x)}}{p(x)} \left( |\nabla U|^{p(x)} + |U|^{p(x)} \right)
- \lambda \int g(a_{i})\frac{s_{0}^{q(x)}}{q(x)}  \left| U \right|^{q(x)}\\
& \qquad - \int f(a_{i}) \frac{s_{0}^{r(x)}}{r(x)} |U|^{r(x)}\\
& \leq  J_{\infty}(s_{0}U) \leq \max_{s \geq 0}J_{\infty} (sU) = J_{\infty}(U) = c_{\infty}.
\end{align*}
Consequently,
\begin{align*}
\limsup_{k \to +\infty}(\sup_{t \geq 0 } J_{\lambda, k}(t\widehat{U}^{i}_{k})) \leq  c_{\infty} \,\,\, \,\, \mbox{for} \,\,\, i \in \{1,....,\ell\}.
\end{align*}

Since $ Q_{k}(\widehat{U}^{i}_{k}) \to a_{i} $ as $ k \to \infty $, then $ \widehat{U}^{i}_{k} \in \theta^{i}_{\lambda, k}  $ for all $ k $ large enough. On the other hand, by Claim \ref{ltda-Jl-coo}, $ J_{\lambda, k} (\widehat{U}^{i}_{k}) < c_{\infty} + \frac{\delta_0}{4} $ holds also for $k$ large enough and $\lambda \in [0, \Lambda^*)$. This way, there exists $k_4 \in \mathbb{N}$ such that
$$
\beta^{i}_{\lambda, k} < c_{\infty} + \frac{\delta_0}{4}, \,\,\,\, \forall \lambda \in [0, \Lambda^*) \,\,\, \mbox{and} \,\,\, k \geq k_4.
$$
Thus, decreasing $\delta_0$ if necessary, we can assume that
$$
\beta^{i}_{\lambda, k} < c_{\infty} + \varrho, \,\,\,\, \forall \lambda \in [0, \Lambda^*) \,\,\, \mbox{and} \,\,\, k \geq k_4.
$$
In order to prove the other inequality, we observe that Lemma \ref{lemK2} yields $ J_{\lambda, k}(u) \geq c_{\infty} + \frac{\delta_{0}}{2} $ for all $ u \in \partial \theta^{i}_{\lambda, k} $, if $\lambda \in [0, \Lambda^*)$ and $k \geq k_3$. Therefore,
\begin{align*}
\tilde{\beta}^{i}_{\lambda, k} \geq c_{\infty} + \frac{\delta_{0}}{2}, \,\,\, \mbox{for} \,\,\, \lambda \in [0, \Lambda^*) \,\,\,\, \mbox{and} \,\,\, k \geq k_3.
\end{align*}
Fixing $k^*=\max\{k_3,k_4\}$, we derive that
\[
\beta^{i}_{\lambda, k} < \tilde{\beta}^{i}_{\lambda, k},
\]
for $ \lambda \in [0, \Lambda_{*}) $ and $ k \geq k^*$. \fim

\begin{lem}\label{PSb}
For each $ 1 \leq i \leq\ell$, there exists a $(PS)_{\beta^{i}_{\lambda, k}}$ sequence,   $ \left\{ u^{i}_{n} \right\} \subset \theta^{i}_{\lambda, k} $ for functional $ J_{\lambda, k} $.
\end{lem}

\noindent \textbf{Proof.} By Lemma~\ref{rho}, we know that $ \beta^{i}_{\lambda, k} < \tilde{\beta}^{i}_{\lambda, k}$. Then, the lemma follows adapting the same ideas explored in \cite{Lin12}. \fim

\section{Proof of Theorem~\ref{T1} }

Let $ \{ u^{i}_{n} \} \subset \theta^{i}_{\lambda, k} $ be a $(PS)_{\beta^{i}_{\lambda, k}} $ sequence   for functional $ J_{\lambda, k} $ given by Lemma~\ref{PSb}. Since $ \beta^{i}_{\lambda, k} < c_{\infty} + \varrho$, by Lemma~\ref{Cond-PS} there is $ u^{i}$ such that $ u^{i}_{n} \to u^{i} $ in $ W^{1, p(x)}(\mathbb{R}^{N}) $. Thus,
$$
u^{i} \in \theta^{i}_{\lambda, k}, \,\,\, J_{\lambda, k}(u^i) = \beta^{i}_{\lambda} \,\,\, \mbox{and} \,\,\, J'_{\lambda, k}(u^i) = 0.
$$
Now, we infer that $ u^{i} \neq u^{j} $ for $ i \neq j $ as $ 1 \leq i,j \leq \ell $. To see why, it remains to observe that
$$
Q_k(u^i) \in \overline{B_{\rho_{0}}(a_{i})} \,\,\, \mbox{and} \,\,\,\, Q_k(u^j) \in \overline{B_{\rho_{0}}(a_{j})}.
$$
Once that
$$
\overline{B_{\rho_{0}}(a_{i})} \cap \overline{B_{\rho_{0}}(a_{j})} = \emptyset \,\,\, \mbox{for} \,\,\, i \not= j,
$$
it follows that $u^i \not= u^j$ for $i \not= j$. From this, $ J_{\lambda, k} $ has at least $ \ell $ nontrivial critical points for $\lambda \in [0, \Lambda^*)$ and $ k \geq k^*$, proving the theorem. \fim

\end{document}